\documentclass[12pt,twoside,leqno]{article}
\usepackage{amsmath}
\usepackage{amssymb}
\usepackage{amsxtra}
\usepackage{amscd}
\usepackage{amsthm}
\usepackage[mathscr]{eucal}

\theoremstyle{plain}

\newtheorem{prop}[subsection]{Proposition}

\newtheorem{lem}[subsection]{Lemma}
\newtheorem{conj}[subsection]{Conjecture}

\theoremstyle{definition}

\newtheorem{para}[subsection]{}

\newenvironment{pf}{\proof[\proofname]}{\endproof}

\begin{document}

\title{Heights of motives}

\author{Kazuya Kato}

\maketitle

\newcommand\Cal{\mathcal}
\newcommand\define{\newcommand}

\define\gp{\mathrm{gp}}%
\define\fs{\mathrm{fs}}%
\define\an{\mathrm{an}}%
\define\mult{\mathrm{mult}}%
\define\Ker{\mathrm{Ker}\,}%
\define\Coker{\mathrm{Coker}\,}%
\define\Hom{\mathrm{Hom}\,}%
\define\Ext{\mathrm{Ext}\,}%
\define\rank{\mathrm{rank}\,}%
\define\gr{\mathrm{gr}}%
\define\cHom{\Cal{Hom}}
\define\cExt{\Cal Ext\,}%

\define\cA{\Cal A}
\define\cC{\Cal C}
\define\cD{\Cal D}
\define\cO{\Cal O}
\define\cS{\Cal S}
\define\cM{\Cal M}
\define\cG{\Cal G}
\define\cH{\Cal H}
\define\cE{\Cal E}
\define\cF{\Cal F}
\define\cN{\Cal N}
\define\fF{\frak F}
\define\Dc{\check{D}}
\define\Ec{\check{E}}

\newcommand{\A}{{\mathbb{A}}}
\newcommand{\N}{{\mathbb{N}}}
\newcommand{\Q}{{\mathbb{Q}}}
\newcommand{\Z}{{\mathbb{Z}}}
\newcommand{\R}{{\mathbb{R}}}
\newcommand{\C}{{\mathbb{C}}}
\newcommand{\bN}{{\mathbb{N}}}
\newcommand{\bQ}{{\mathbb{Q}}}
\newcommand{\bF}{{\mathbb{F}}}
\newcommand{\bZ}{{\mathbb{Z}}}
\newcommand{\bP}{{\mathbb{P}}}
\newcommand{\bR}{{\mathbb{R}}}
\newcommand{\bC}{{\mathbb{C}}}
\newcommand{\bbQ}{{\bar \mathbb{Q}}}
\newcommand{\ol}[1]{\overline{#1}}
\newcommand{\too}{\longrightarrow}
\newcommand{\respect}{\rightsquigarrow}
\newcommand{\compatible}{\leftrightsquigarrow}
\newcommand{\upc}[1]{\overset {\lower 0.3ex \hbox{${\;}_{\circ}$}}{#1}}
\newcommand{\Gmlog}{\bG_{m, \log}}
\newcommand{\Gm}{\bG_m}
\newcommand{\ep}{\varepsilon}
\newcommand{\Spec}{\operatorname{Spec}}
\newcommand{\val}{{\mathrm{val}}} 
\newcommand{\n}{\operatorname{naive}}
\newcommand{\bs}{\operatorname{\backslash}}
\newcommand{\Gal}{\operatorname{{Gal}}}
\newcommand{\gal}{{\rm {Gal}}({\bar \Q}/{\Q})}
\newcommand{\galp}{{\rm {Gal}}({\bar \Q}_p/{\Q}_p)}
\newcommand{\gall}{{\rm{Gal}}({\bar \Q}_\ell/\Q_\ell)}
\newcommand{\wep}{W({\bar \Q}_p/\Q_p)}
\newcommand{\wel}{W({\bar \Q}_\ell/\Q_\ell)}
\newcommand{\Ad}{{\rm{Ad}}}
\newcommand{\BS}{{\rm {BS}}}
\newcommand{\even}{\operatorname{even}}
\newcommand{\End}{{\rm {End}}}
\newcommand{\odd}{\operatorname{odd}}
\newcommand{\GL}{\operatorname{GL}}
\newcommand{\np}{\text{non-$p$}}
\newcommand{\g}{{\gamma}}
\newcommand{\G}{{\Gamma}}
\newcommand{\Lam}{{\Lambda}}
\newcommand{\La}{{\Lambda}}
\newcommand{\lam}{{\lambda}}
\newcommand{\la}{{\lambda}}
\newcommand{\uL}{{{\hat {L}}^{\rm {ur}}}}
\newcommand{\uQp}{{{\hat \Q}_p}^{\text{ur}}}
\newcommand{\sel}{\operatorname{Sel}}
\newcommand{\dt}{{\rm{Det}}}
\newcommand{\Sig}{\Sigma}
\newcommand{\fil}{{\rm{fil}}}
\newcommand{\SL}{{\rm{SL}}}
\newcommand{\spl}{{\rm{spl}}}
\newcommand{\st}{{\rm{st}}}
\newcommand{\Isom}{{\rm {Isom}}}
\newcommand{\Mor}{{\rm {Mor}}}
\newcommand{\bg}{\bar{g}}
\newcommand{\id}{{\rm {id}}}
\newcommand{\cone}{{\rm {cone}}}
\newcommand{\al}{a}
\newcommand{\ChL}{{\cal{C}}(\La)}
\newcommand{\Image}{{\rm {Image}}}
\newcommand{\toric}{{\operatorname{toric}}}
\newcommand{\torus}{{\operatorname{torus}}}
\newcommand{\Aut}{{\rm {Aut}}}
\newcommand{\Qp}{{\mathbb{Q}}_p}
\newcommand{\barQp}{{\mathbb{Q}}_p}
\newcommand{\Qpur}{{\mathbb{Q}}_p^{\rm {ur}}}
\newcommand{\Zp}{{\mathbb{Z}}_p}
\newcommand{\Zl}{{\mathbb{Z}}_l}
\newcommand{\Ql}{{\mathbb{Q}}_l}
\newcommand{\Qlur}{{\mathbb{Q}}_l^{\rm {ur}}}
\newcommand{\F}{{\mathbb{F}}}
\newcommand{\eps}{{\epsilon}}
\newcommand{\epsLa}{{\epsilon}_{\La}}
\newcommand{\epsLaVxi}{{\epsilon}_{\La}(V, \xi)}
\newcommand{\epsOLaVxi}{{\epsilon}_{0,\La}(V, \xi)}
\newcommand{\Qplin}{{\mathbb{Q}}_p(\mu_{l^{\infty}})}
\newcommand{\otimesQplin}{\otimes_{\Qp}{\mathbb{Q}}_p(\mu_{l^{\infty}})}
\newcommand{\galFl}{{\rm{Gal}}({\bar {\Bbb F}}_\ell/{\Bbb F}_\ell)}
\newcommand{\gallur}{{\rm{Gal}}({\bar \Q}_\ell/\Q_\ell^{\rm {ur}})}
\newcommand{\galFF}{{\rm {Gal}}(F_{\infty}/F)}
\newcommand{\galFv}{{\rm {Gal}}(\bar{F}_v/F_v)}
\newcommand{\galF}{{\rm {Gal}}(\bar{F}/F)}
\newcommand{\epsVxi}{{\epsilon}(V, \xi)}
\newcommand{\epsOVxi}{{\epsilon}_0(V, \xi)}
\newcommand{\plim}{\lim_
{\scriptstyle 
\longleftarrow \atop \scriptstyle n}}
\newcommand{\sig}{{\sigma}}
\newcommand{\ga}{{\gamma}}
\newcommand{\del}{{\delta}}
\newcommand{\Vss}{V^{\rm {ss}}}
\newcommand{\Bst}{B_{\rm {st}}}
\newcommand{\Dpst}{D_{\rm {pst}}}
\newcommand{\Dcrys}{D_{\rm {crys}}}
\newcommand{\DdR}{D_{\rm {dR}}}
\newcommand{\Fin}{F_{\infty}}
\newcommand{\Kla}{K_{\lambda}}
\newcommand{\Ola}{O_{\lambda}}
\newcommand{\Mla}{M_{\lambda}}
\newcommand{\Det}{{\rm{Det}}}
\newcommand{\Sym}{{\rm{Sym}}}
\newcommand{\LaSa}{{\La_{S^*}}}
\newcommand{\cX}{{\cal {X}}}
\newcommand{\MHG}{{\frak {M}}_H(G)}
\newcommand{\tauMla}{\tau(M_{\lambda})}
\newcommand{\Fvur}{{F_v^{\rm {ur}}}}
\newcommand{\Lie}{{\rm {Lie}}}
\newcommand{\cL}{{\cal {L}}}
\newcommand{\cW}{{\cal {W}}}
\newcommand{\fq}{{\frak {q}}}
\newcommand{\cont}{{\rm {cont}}}
\newcommand{\SC}{{SC}}
\newcommand{\Om}{{\Omega}}
\newcommand{\dR}{{\rm {dR}}}
\newcommand{\crys}{{\rm {crys}}}
\newcommand{\hatSig}{{\hat{\Sigma}}}
\newcommand{\rdet}{{{\rm {det}}}}
\newcommand{\ord}{{{\rm {ord}}}}
\newcommand{\BdR}{{B_{\rm {dR}}}}
\newcommand{\BdRO}{{B^0_{\rm {dR}}}}
\newcommand{\Bcrys}{{B_{\rm {crys}}}}
\newcommand{\Qw}{{\mathbb{Q}}_w}
\newcommand{\barkappa}{{\bar{\kappa}}}
\newcommand{\cP}{{\Cal {P}}}
\newcommand{\cZ}{{\Cal {Z}}}
\newcommand{\oppLa}{{\Lambda^{\circ}}}

\begin{abstract}
 We define the height of a motive over a number field. We show that if we assume  the finiteness of motives of bounded height, Tate conjecture for the $p$-adic Tate module  can be proved for motives with good reduction at $p$. 
\end{abstract}
\renewcommand{\thefootnote}{\fnsymbol{footnote}}
\footnote[0]{The author is partially supported by an NSF grant.}

\begin{para}\label{0.1}

In this paper, we generalize the definition of the height of an abelian variety over number field due to Faltings \cite{Fa} to a motive over a number field. Here by motive, we mean a pure motive.

We define the height of a motive $M$ over a number field $K$  as the 
Arakelov degree of the one dimensional $\Q$-vector space $$L(M)_{\Q}:= \otimes_{r\in \Z}\;   (\text{det}_\Q\; \gr^r M_{dR})^{\otimes r}$$
which is endowed with a metric and an integral structure. Here $M_{dR}$ denotes the de Rham realization of $M$, and $\gr^r$ is that of  the Hodge filtration on $M_{dR}$.

The metric on $L(M)_{\Q}$ is defined by using Hodge theory (see section \ref{Ho}), and the integral structure of $L(M)_{\Q}$ is defined by using $p$-adic Hodge theory (see section \ref{pHo}).

In the case $M$ is the $H_1$ of an abelian variety $A$, $\gr^rM_{dR}$ is $Lie(A)$ if $r=-1$ and is $0$ if $r\neq 0, -1$, and our definition coincides with the height of Faltings who used the Arakelov degree of  $coLie(A)$ (the dual of $Lie(A)$).

Height is a basic notion in number theory. We hope that our generalization of this notion to motives will supply
fruitful subjects and interesting problems to number theory.

Details of this paper will be given elsewhere. 

\end{para}

\section{Hodge theory}\label{Ho}

\begin{para}

In this section, for a pure Hodge structure $H= (H_{\Z}, F)$, we define a canonical metric on the one dimensional $\C$-vector space $$L(H):=\otimes_{r\in \Z}\; ({\det}_\C\;  F^r/F^{r+1})^{\otimes r}.$$

\end{para}
\begin{para}
Let $H=(H_{\Z}, F)$ be a pure Hodge structure of weight $w$. That is, $H_\Z$ is a free $\Z$-module of finite rank and $F$ is a descending filtration on $H_\C=\C\otimes H_\Z$ such that $F^r=H_\C$ for $r\ll 0$ and $F^r=0$ for $r\gg 0$, satisfying 
$$H_\C=\oplus_{r\in \Z} \; H_{\C}^{r,w-r}\quad \text{where}\;\; H_\C^{r,w-r}= F^r \cap {\bar F}^{w-r}.$$
Here ${\bar F}^r$ is the image of $F^r$ under the complex conjugate $H_\C\to H_\C\;;\;z\otimes h\mapsto {\bar z}\otimes h$ ($z\in \C, h\in H_\Z$). 

Let 
$$L'(H):=\otimes_{r\in \Z}\;({\det}_\C\; H^{r,w-r})^{\otimes r}.$$ The canonical isomorphism $H^{r,w-r}\overset{\cong}\to F^r/F^{r+1}$ induces a canonical isomorphism $L'(H)\overset{\cong}\to L(H)$.

\end{para}

\begin{para}

The complex conjugate $H_\C\to H_\C$ induces an isomorphism $$L'(H) \overset{\cong}\to \bar L'(H):= \otimes_{r\in \Z}\;({\det}_\C\;  H^{w-r,r})^{\otimes r}=\otimes_{r\in \Z}\;({\det}_\C\;  H^{r,w-r})^{\otimes w-r}$$ where the last $=$ is obtained by replacing $r$ by $w-r$. We have a canonical isomorphism
$$L'(H)\otimes_{\C} \bar L'(H)\cong \C\otimes_\Z ({\det}_\Z\;  H_\Z)^{\otimes w}$$
as follows:
$$L'(H)\otimes_{\C} \bar L'(H)= (\otimes_{r\in \Z}\; ({\det}_\C\;  H^{r,w-r})^{\otimes r}) \otimes (\otimes_{r\in \Z}\; ({\det}_\C\;  H^{r,w-r})^{\otimes w-r})$$ $$
= \otimes_{r\in \Z} \; ({\det}_\C\;  H^{r,w-r})^{\otimes w}= ({\det}_\C\;  H_\C)^{\otimes w}=\C\otimes_\Z ({\det}_\Z\;  H_\Z)^{\otimes w}.$$

\end{para}

\begin{para}

Let $s\in L'(H)$. Then we have an element $\bar s$ of $\bar L'(H)$, and $s\otimes \bar s$ is sent to an element $ze$ of $\C\otimes_\Z\; ({\det}_\Z\;  H_\Z)^{\otimes w}$ via the canonical isomorphism, where $z\in \C$ and $e$ is a $\Z$-basis of $({\det}_\Z\;  H_\Z)^{\otimes w}$. We define  $|s|=|z|^{1/2}$. 

Via the canonical isomorphism $L'(H)\overset{\cong}\to L(H)$, we obtained a metric on $L(H)$. 
\end{para}

\section{$p$-adic Hodge theory}\label{pHo}

\begin{para} In this section, let $p$ be a prime number and let $K$ be a finite extension of $\Q_p$. Let $T$ be a free $\Z_p$-module  of finite rank endowed with a continuous action of $G_K:=\Gal(\bar K/K)$ such that $V= \Q_p\otimes_{\Z_p} T$ is de Rham in the sense of Fontaine \cite{Fo}. The goal of this section is to define a $p$-adic integral structure (a $\Z_p$-structure) $L_r(T)$ of $$L_r(V):= {\det} _{\Q_p} (D_{dR}(V)/D^r_{dR}(V))$$ for each $r\in \Z$. 

$L_r(T)$ is defined to be $L_0$ of the Tate twist $T(r)$. So, we consider $L_0(T)$.
\end{para}

\begin{para}  
As in \cite{BK}, we have an exact sequence 
$$0\to H^0(K, V) \to D_{\crys}(V) \to D_{\crys}(V) \oplus D_{dR}(V)/D_{dR}^0(V) \to H^1_f(K, V) \to 0$$
where $H^m$ are Galois cohomology and $H^1_f(K, V)\subset H^1(K, V)$ is as in \cite{BK}. The map 
$D_{crys}(V)\to D_{crys}(V)$ is $x\mapsto (1-\varphi)x$ where $\varphi$ is the Frobenius, and the map $D_{\crys}(V)\to D_{dR}(V)/D_{dR}^0(V)$ is the evident map. 

We will define a $\Z_p$-submodule $H_{cf}^1(K, T)$ of $H^1_f(K, T)$ of finite index, and define the $p$-adic integral structure $L_0(T)$ of $L_0(V)$ as
$$L_0(T):= \text{det}_{\Z_p} H^1_{cf}(K, T)\otimes (\text{det}_{\Z_p}H^0(K, T))^{\otimes -1}
\subset {\det}_{\Q_p} \; H^1_f(K, V)\otimes ({\det}_{\Q_p}\; H^0(K, V))^{\otimes -1}
$$ $$\cong  {\det}_{\Q_p}\; D_{dR}(V)/D^0_{dR}(V) \otimes {\det}_{\Q_p}\; D_{\crys}(V) \otimes ({\det}_{\Q_p}\; D_{\crys}(V))^{\otimes -1}
\cong {\det}_{\Q_p}\;  D_{dR}(V)/D^0_{dR}(V)$$
where the first isomorphism is obtained by the above exact sequence and the second isomorphism is obtained by canceling two $D_{\crys}(V)$ via the identity map.

\end{para}

\begin{para} We describe our motivation of the definition of $H^1_{cf}(K, T)$.

Let $A$ be an abelian variety over $K$, and let $TA=\prod_{\ell} T_{\ell}A$ where $\ell$ ranges over all prime numbers (including $p$) and $T_{\ell}A$ is the $\ell$-adic Tate module. Then the Kummer sequences 
$$0\to TA/nTA \to A(\bar K) \overset{n}\to A(\bar K) \to 0$$
for $n\geq 1$ induce  connecting  homomorphisms $A(K) \to H^1(K, TA/nTA)$, and the inverse limit gives an isomorphism $A(K)\overset{\cong}\to H^1_f(K, TA)$ (see \cite{BK}). 

Let $\cA$ be the N\'eron model of $A$ and let $\cA^{\circ}\subset \cA$ be the connected N\'eron model of $A$ ($\cA^{\circ}$  is the open set of $\cA$ obtained from $\cA$ by removing all connected components of the special fiber of $\cA$ which do not contain the origin). 
Thus $$\cA^{\circ}(O_K) \subset \cA(O_K)=A(K)\cong H^1_f(K, TA).$$
For our seek of the nice integral structure on the de Rham object, $\cA^{\circ}(O_K)$ is important. We can identify  $\cA^{\circ}(O_K)$  as the subgroup $H_{cf}^1(K, TA)$ of $H^1_f(K, TA)$.

\end{para}

\begin{prop}\label{prop}  Let $\ell$ be a prime number, and let $T$ be a free $\Z_{\ell}$-module of finite rank endowed with a continuous action of $\Gal(\bar K/K)$. In the case $\ell=p$, assume that $V=\Q_{\ell}\otimes_{\Z_{\ell}} T$ is de Rham. Let $a\in H^1_f(K,T)$. For $n\geq 1$, let $K_n$ be the unique unramified extension of $K$ of degree $n$.

\medskip

(1) The following three conditions are equivalent.

\medskip

(a) For any $n\geq 1$, $a$ belongs to the image of the trace map $H^1_f(K_n, T)\to H_f^1(K, T)$. 

\smallskip

(b)  For any $n\geq 1$, the image of $a$ in $H^1_f(K_{ur},T)$ belongs to the image of $1-\varphi^n: H^1_f(K_{ur}, T) \to H^1_f(K_{ur}, T)$. Here $K_{ur}$ is the maximal unramified extension of $K$ and $\varphi$ is
the Frobenius of $K_{ur}/K$.

\smallskip

(c) For any $n$, the map $H^1_f(K_n, T_a)\to H^1_f(K_n, \Z_{\ell})=\Z_{\ell}$ is surjective. Here $T_a$ is the extension of $\Z_p$ by $T$ corresponding to $a$. 

\medskip

(2) If $\ell\neq p$, the equivalent conditions (a) -- (c) are also equivalent to
\medskip

(d)  $a$ belongs to the kernel of $H^1(K, T)\to H^1(K_{ur}, T)$. 

\medskip

(3) If $T=T_{\ell}A$ for an abelian vareity $A$ over $K$, the equivalent conditions (a)--(c) are also equivalent to

\medskip

(e) $a$ belongs to the image of $\cA^{\circ}(O_K)$. 
\end{prop}

Concerning (2), note that in the case $\ell\neq p$, $H^1_f(K, T)$ is defined to be the kernel of $H^1(K, T)\to H^1(K_{ur}, V)$  which can be a little bigger than the kernel of $H^1(K, T)\to H^1(K_{ur}, T)$.

\begin{para} We define 
$H^1_{cf}(K, T)$ to be the subgroup of $H^1_f(K, T)$ consisting of all elements $a$ satisfying the equivalent conditions in \ref{prop} (1).  We call it the connected finite part of $H^1(K, T)$. 
We can prove that $H^1_f(K, T)/H^1_{cf}(K, T)$ is finite.
\end{para}

\begin{para} Let $A$ be an abelian variety over $K$ and let $T=T_pA$. Then we have $$L_0(T)=\text{det}_{\Z_p}\; \text{Lie}(\cA)= \text{det}_{\Z_p}\;\text{Lie}(\cA^{\circ}).$$

\end{para}

\begin{para}\label{goodp} Assume $K$ is unramifed over $\Q_p$, $V$ is crystalline, and there is $a\in \Z$ such that $D^a_{dR}(V)=D_{dR}(V)$ and $D^{a+p-1}_{dR}=0$. Then there is an $O_K$-lattice $D$ of $D_{dR}(V)=D_{crys}(V)$  corresponding to $T$ which satisfies  $D=\sum_{i\in \Z}\; p^{-i}\varphi D^i$ where $D^i=D\cap D^i_{dR}(V)$. 
($T$ is constructed from $D$ by the method of Fontaine-Laffaile.) In this case, we have
$$H^1_f(K, T)=H^1_{cf}(K, T)$$
and we have an exact sequence
$$0\to H^0(K, T) \to D^0\overset{1-\varphi}\to D\to H^1_{cf}(K, T)\to 0.$$
In the case $T=H^m_{et}(X_{\bar K}, \Z_p)(r)/({\text{torsion}})$ for a proper smooth scheme $X$ over $O_K$ with $m\in \Z$  such that $m\leq p-2$ and with $r\in \Z$, if $Y$ denotes the special fiber of $X$, then
$D=H^m_{dR}(X/O_K)/({\text{torsion}})=H^m_{\crys}(Y)/({\text{torsion}})$, $(D^i)_i$ coincides with $r$-twist of the Hodge filtration of $H^m_{dR}(X/O_K)/(\text{torsion})$, and $\varphi$ of $D$ coincides with $p^{-r}\varphi$ of $H^m_{\crys}(Y)/(\text{torsion})$. 

\end{para}

\begin{para} Let $K'$ be a finite extension of $\Q_p$ contained in $K$ and let $T'$ be the Weil restriction of $T$ to $K'$ (that is, the induced representation of $G_{K'}$ obtained from $T$). Then $H^1_f(K, T)=H^1_f(K', T')$, $H^1_{cf}(K, T)=H^1_{cf}(K', T')$, and $L_0(T)=L_0(T')$. 

This corresponds to the fact that the Weil restriction of $\cA$ (resp. $\cA^{\circ}$) to $O_{K'}$ is the N\'eron model (resp. connected N\'eron model) of the Weil restriction $A'$ of $A$ to $K'$. 
\end{para}

\section{Heights of motives}

\begin{para} 

Let $\A^f$ be the ring of finite adeles of $\Q$.

Let $k$ be a field of characteristic $0$. For a motive $M$ over $k$ (this means the usual pure motive with $\Q$-coefficients), let $M_{\A^f}$ is the \'etale realization of $M$ with $\A^f$-coefficients which is endowed with the continuous action of $G_k=\Gal(\bar k/k)$.

By a $\Z$-motive over $k$, we mean a motive $M$ over $k$ endowed with a  $G_k$-stable  $\hat \Z$-submodule $T$ of $M_{\A^f}$ which is free of finite type over $\hat \Z$ such that $\A^f\otimes_{\hat \Z} T = M_{\A^f}$. 

In this section, we define the height of a $\Z$-motive over a number field. 

\end{para}

\begin{para}
To avoid a technical problem, we fix integers $a, b$ such that $a\leq b$, and we define the height of a $\Z$-motive $M$ over a number field $K$ satisfying $M^a_{dR}=M_{dR}$ and $M^b_{dR}=0$, depending on the choices of $a, b$. 

We define the height of such $M$ as the height of the Weil restriction of $M$ to $\Q$. Hence we consider $\Z$-motives over $\Q$. 
\end{para}

\begin{para}
Let $M=(M, T)$ be a $\Z$-motive over $\Q$ such that $M_{dR}^a=M_{dR}$ and $M_{dR}^b=0$. Let
$$L(M)_{\Q} = \otimes_{r\in \Z} \; (\text{det}_{\Q} \; \gr^r M_{dR})^{\otimes r}.$$

We define a  metric on $L(M)_{\R}=\R\otimes_{\Q} L(M)_\Q$. Let $M_B$ be the Betti realization which is a $\Q$-vector space.  Then we have a $\Z$-structure $M_{B,\Z}$ of $M_B$ by
$$M_{B, \Z} = M_B \cap T \subset M_{\A^f}.$$ Via the canonical ismorphism $M_{B, \C}\cong M_{dR, \C}$, $M_{B, \C}=\C\otimes_{\Z} M_{B, \Z}$ has a pure Hodge filtration
$F^r:=\C\otimes_{\Q} M^r_{dR}$ ($r\in \Z$), and $H=(M_{B, \Z}, F)$ is a Hodge structure. Hence by section \ref{Ho}, we have a metric on $L(H)=L(M)_{\C}$. By restricting to $L(M)_{\R}$. we have a metric on $L(M)_{\R}$.

Next we define an integral structure $L(M)_{\Z}$ on $L(M)_{\Q}$. 
For each $r\in \Z$, let
$$L_r(M)_{\Q} = \text{det}_{\Q} (M_{dR}/M^r_{dR}).$$
 Then by section 2, the $p$-adic component $T_p$ of $T$ regarded as a representation of $G_{\Q_p}$ defines a $\Z_p$-structure $L_r(T_p)$ of $L_r(M)_{\Q_p}$. When $p$ ranges, this gives a $\hat \Z$-structure $L_r(M)_{\hat \Z}$ of $\A^f\otimes_{\Q} L(M)_{\Q}$. This follows from the fact that if $(M, T)= (H^m(X), H^m_{et}(X_{\bar \Q}, \hat \Z))$ for a projective smooth scheme $X$ over $\Q$ and if $\frak X$ is a projective scheme over $\Z$ such that $X=\frak X \otimes \Q$, then $L_r(T_p)=\Z_p\otimes_{\Z}  {\det}_{\Z} (H^m_{dR}(\frak X/\Z)/F^rH^m_{dR}(\frak X/\Z))$ for almost all $p$ by \ref{goodp}.
 Let
$$L_r(M)_{\Z} = L_r(M)_{\Q}\cap L_r(M)_{\hat \Z} \subset \A^f\otimes_{\Q} L_r(M)_\Q.$$

We define
$$L(M)_{\Z} := (\otimes_{a<i<b} L_i(M)_{\Z}^{\otimes -1})\otimes L_b(M)_{\Z}^{\otimes (b-1)}\subset
(\otimes_{a<i<b} L_i(M)_{\Q}^{\otimes -1})\otimes L_b(M)_{\Q}^{\otimes (b-1)}=L(M)_{\Q}.$$

The reason why we do not take the simpler definition $$L(M)_{\Z} := \otimes_{r\in \Z} \; (L_r(M)_{\Z}^{\otimes -1}\otimes L_{r+1}(M)_{\Z})^{\otimes r}$$ (independently of $a$ and $b$) is that we are not sure whether this is a finite tensor product. 
\end{para}
\begin{para}\label{h(M)}
For a $\Z$-motive $M$ over $\Q$, we define its multiplicative height $H(M)$ and the logarithmic height $h(M)$ as
$$H(M) = |e|^{-1}, \quad h(M) =- \log (|e|)$$
where $e$ is a $\Z$-basis of $L(M)_{\Z}$ and $|\;|$ is the metric of $L(M)_{\R}$.
\end{para}

\begin{para}
If $A$ is an abelian variety over a number field $K$, for $M=(H_1(A), T(A))$, the height of $A$ defined by Faltings coincides with our height for the choice $a=-1$ and $b=1$. 
\end{para}

\section{Some topics}\label{consec}

\begin{para} 
We fix the type $\Phi=(w, (h^r)_{r\in \Z})$ of motives, where $w$ is the weight and $h^r=\dim \gr_{dR}^r$.  Take $a, b$ such that $h^r=0$ unless $a\leq r < b$ and define the height of a $\Z$-motive by using these fixed $a,b$.

\end{para}
The following is a basic conjecture.
\begin{conj}\label{conj1}   Let $K$ be a number field and let $c>0$. Then there are only finitely many isomorphism classes of motives over $K$ of type $\Phi$ of semi-stable reduction such that $h(M)\leq c$.
\end{conj}

\begin{para}

In \cite{Fa}, by using his heights of abelian varieties, Faltings proved the Tate conjecture $\Z_p \otimes Hom(A, B) \overset{\cong}\to  Hom_{G_K}(T_p A, T_p B)$ for abelian varieties $A$ and $B$ over a number field $K$. A key point was that the above finiteness is true for abelian varieties. He proved this finiteness by using the fact that his height of an abelian variety essentially coincides with the height of the corresponding point of the moduli space of abelian varieties. For the above general conjecture, the difficulty is that usually there is no moduli space of $\Z$-motives of type $\Phi$. 
\end{para}

\begin{prop}\label{Tate} Let $M=(M, T)$ and $M'=(M', T')$ be $\Z$-motives over a number field $K$ of type $\Phi$. Let $p$ be a prime number,  and assume that 
$M_p$ and $M'_p$ are crystalline as representations of $G_{K_v}$ for any place $v$ of $K$ lying over $p$, $p$ is unramified in $K/\Q$, and  $b\leq a+p-1$. 
Assume that Conjecture \ref{conj1} is true. Then
$$\Z_p\otimes_{\Z} Hom(M, M') \overset{\cong}\to Hom_{G_K}(T_p, T'_p).$$
\end{prop}

 The is proved in the following way. 

\begin{lem}\label{lem} Let $V$ be a finite dimensional $\Q_p$-vector space endowed with a continuous action of $G_\Q$. Assume that $V$ is de Rham as a representation of $G_{\Q_p}$ and assume that there is an integer $w\in \Z$ such that for almost all prime numbers $\ell$, the action of $G_{\Q_{\ell}}$ on $V$ is unramified and all eigen values of the action of the geometric frobenius of $\ell$ on $V$ are algebraic numbers whose all conjugates are of complex absolute value $\ell^{w/2}$. Let
$$s(V) = \sum_{r\in \Z}\; r\cdot \dim \gr^rD_{dR}(V),   \quad t(V)=w\cdot \dim(V)/2.$$ Then we have $s(V)=t(V)$.
\end{lem}
\begin{pf} This is reduced to the case $\dim(V)=1$ by the facts $s(V)=s(\det(V))$ and $t(V)=t(\det(V))$. If $\dim(V)=1$,  $V$ is isomorphic to $\Q_p(m)$ for some integer $m$ as a representation of $G_K$ for some finite extension $K$ of $\Q$, and hence $s(V)=-m=t(V)$.
\end{pf}

If we assume Conjecture \ref{conj1}, 
then by the argument of Faltings in \cite{Fa}, the proof of Proposition \ref{Tate} is  reduced to 
\begin{prop} Let $M=(M, T)$ be as in the hypothesis of Proposition \ref{Tate}. (We do not assume Conjecture \ref{conj1} here.) Let $U$ be a free $\Z_p$-module of finite rank endowed with an action of $G_K$, and assume that we have a surjective homomorphism $T_p\to U$ which is compatible with the actions of $G_K$. For $n\geq 0$, let $T^{(n)}:= Ker(T\to U/p^nU)$, and let $M^{(n)}$ be the $\Z$-motive over $K$ which is the same as $M$ as a $\Q$-motive over $K$ but with the Galois representation $T^{(n)}$ over $\hat \Z$. 
Then
$h(M^{(n)})=h(M)$ for any $n\geq 0$. 
\end{prop}

\begin{pf}  By Weil restriction, we may assume $K=\Q$. Let $D$ be as in  \ref{goodp}. 
We have exact sequences $0\to D^i(T_p^{(n)})\to D^i(T_p) \to D^i(U)/p^nD^i(U)\to 0$ for all $i$.  From this, we have
$$L(M^{(n)})_{\Z}=p^{ns(V)}L(M)_{\Z}$$ where $V=\Q_p\otimes_{\Z_p} U$. 
On the other hand, $$({\det}_{\Z} \; H^{(n)}_\Z)^{\otimes w}= p^{2nt(V)}({\det}_{\Z}\; H_\Z)^{\otimes w}.$$
where $H$ (resp. $H^{(n)}$) is the Hodge structure associated to $M$ (resp. $M^{(n)}$). Hence 
$$H(M)/H(M^{(n)}) = p^{n(s(V)-t(V))}=1$$
by Lemma \ref{lem}. 
\end{pf}

\begin{para}
Many questions arise concerning  heights of motives. For example, we have the following analogue of abc conjecture (or Vojta conjecture \cite{Vo}) for motives. For a motive $M$ over a number field $K$, let $n(M)=\sum_v \; \log N(v)$ where $v$ ranges over all finite places of $K$ at which $M$ has bad reduction.

\end{para}
\begin{conj}\label{abc}  There are constants $c, c'>0$ such that $$n(M)\geq c\cdot h(M)- \log(|D_K|) - c'\cdot [K:\Q]$$ for any number field $K$ and any $\Z$-motive $M$ over $K$ of type $\Phi$ of semi-stable reduction. Here $D_K$ denotes the discriminant of $K$. 
\end{conj}

Kazuya KATO,

Department of mathematics,
University of Chicago,
Chicago, Illinois, 60637,
USA,
{\tt kkato@math.uchicago.edu}

\end{document}